\documentclass{amsart}
\usepackage{amsmath}
\usepackage{amsthm}
\usepackage{amsfonts}
\usepackage{amssymb}
\usepackage{amscd}
\usepackage{cite}

\usepackage{graphicx}

\theoremstyle{plain}

\newtheorem*{alt-thm}{Theorem}
\newtheorem*{alt-cor}{Theorem}

\theoremstyle{definition}

\begin{document}

\title[Unit-Regularity of Regular Nilpotent Elements]{Unit-Regularity of Regular Nilpotent Elements}

\author{Dinesh Khurana}
\address{Department of Mathematics, Panjab University, Chandigarh-160014, India}
\email{dkhurana@pu.ac.in}


\begin{abstract}
 Let $a$ be a regular element of a ring $R$. If either $K:=\rm{r}_R(a)$ has the exchange property or every power of $a$ is regular, then we prove that for every positive integer $n$ there exist decompositions $$ R_R = K \oplus X_n \oplus Y_n = E_n \oplus X_n \oplus aY_n,$$ where $Y_n \subseteq a^nR$ and $E_n \cong R/aR$. As applications we get easier proofs of the results that a strongly $\pi$-regular ring has stable range one and also that a strongly $\pi$-regular element whose every power is regular is unit-regular.
\end{abstract}

\maketitle

An element $a$ of a ring $R$ is called strongly $\pi$-regular if both chains $aR \supseteq a^2R \supseteq a^3R...$ and $Ra \supseteq Ra^2 \supseteq Ra^3...$ stabilize. If every element of $R$ is strongly $\pi$-regular, then $R$ is called a {\em strongly $\pi$-regular} ring. In \cite{Ara} Pere Ara proved a wonderful result that a strongly $\pi$-regular ring has stable range one. Ara's proof is on the following lines. As a strongly $\pi$-regular ring is an exchange ring and an exchange ring has stable range one if and only if every regular element is unit-regular, it is enough to show that every regular element of a strongly $\pi$-regular ring is unit-regular. Suppose $a$ is a regular element of a strongly $\pi$-regular ring. By \cite[Proposition 1]{Ni} there exist $n \in \mathbb{N}$, an idempotent $e$ and a unit $u$ in $R$ with $a^n = eu$ such that $a,\,e$ and $u$ commute with each other. Then $ea$ is a unit in $eRe$ with inverse $ea^{n-1}u^{-1}$ and $(1-e)a$ is a regular nilpotent element of the exchange ring $(1-e)R(1-e)$. As $a = ea + (1-e)a$ and $ea$ is unit-regular in $eRe$, we will get that $a$ is unit-regular if we can show that $(1-e)a$ is unit-regular in $(1-e)R(1-e)$. So the result will follow if we can show that a regular nilpotent element of an exchange ring  is unit-regular. This is the crucial result proved by Ara in \cite{Ara} and an easier proof of this will follow from our Theorem 2.

\bigskip
In \cite[Theorem 5.8]{GM} Goodearl and Menal proved that a regular strongly $\pi$-regular ring is unit-regular. The proof of \cite[Theorem 5.8]{GM} can be adapted to prove that if $a$ is a strongly $\pi$-regular element of any ring $R$ such that $a^n$ is regular for every $n \in \mathbb{N}$, then $a$ is unit-regular. A different proof of this result was given by Beidar, O'Meara and Raphael in \cite[Corollary 3.7]{BMR}.  Suppose $a$ is a strongly $\pi$-regular element of a ring $R$ whose each power is regular. Then as above there exist $n \in \mathbb{N}$, an idempotent $e$ and a unit $u$ in $R$ with $a^n = eu$ such that $a,\,e$ and $u$ commute with each other. As seen above it will follow that $a$ is unit-regular if we can prove that $(1-e)a$ is unit-regular in $(1-e)R(1-e)$. As $(1-e)a$ is nilpotent and its each power is regular in $(1-e)R(1-e)$, it is enough to prove that a nilpotent element whose each power is regular is unit-regular. An easier proof of this will follow from our Theorem 4.

\bigskip
Recently Ara and O'Meara in \cite{AO} and Pace and \v{S}ter in \cite{NS} have shown that a regular nilpotent element in general may not be unit-regular.

\bigskip
By $A \subseteq^{\oplus} B$ we shall mean that $A$ is a summand of the module $B$. We will tacitly use the fact that a regular element $a\in R$ is unit-regular if and only if $\rm{r}_R(a) \cong R/aR$, where $\rm{r}_R(a)=\{x\in R: ax=0\}$.

\bigskip
\noindent{\bf Lemma 1 \cite[Corollary 3.9]{CJ}.} {\em If $M$ has the exchange property and
$A=M\oplus B \oplus C = \bigoplus_{I}A_i \oplus C$, then there exists a decomposition $A_i = D_i \oplus E_i$ of each $A_i$ such that $A = M \oplus \bigoplus_{I}D_i \oplus C.$}

\bigskip
\noindent{\bf Theorem 2.} {\em Let $a$ be a regular element of a ring $R$ such that the right $R$-module $K := \rm{r}_R(a)$ has the exchange property. Then for every $n \in \mathbb{N}$ there exist decompositions $$ R_R = K \oplus X_n \oplus Y_n = E_n \oplus X_n \oplus aY_n,$$ where $Y_n \subseteq a^nR$ and $E_n \cong R/aR$. If $a^n = 0$, then $aY_n = Y_n = 0$ and so $K \cong E_n \cong R/aR$ implying that $a$ is unit-regular.}\\[3mm]
{\bf Proof.} For every positive integer $i$ we will inductively construct right ideals $A_i,\;A_i',\;Y_i$ of $R$ such that for every $j\geq 1$, $$R = K \oplus (\bigoplus_{i=1}^{j}A_i) \oplus Y_j = (\bigoplus_{i=1}^{j}A_i) \oplus (A_j' \oplus aA_j) \oplus aY_j, \hspace{10mm} (*)$$ where $Y_j \subseteq a^jR$ and $A_j' \oplus aA_j = A_{j+1} \oplus A_{j+1}' \cong R/aR$. Then we have the desired decompositions by putting $X_n = \bigoplus_{i=1}^{n}A_i$ and $E_n = A_n' \oplus aA_n$.

\bigskip
We have $R = K \oplus B = A \oplus aR$ for some right ideals  $A$ and $B$ of $R$. As $K_R$ has the exchange property, by Lemma 1 we have decompositions $A = A_1 \oplus A_1'$ and $aR = Y_1 \oplus Y_1'$ such that $R = K \oplus A_1 \oplus Y_1$. As $K \cap (A_1 \oplus Y_1) = 0$, $aR = aA_1 \oplus aY_1$ and $aA_1 \cong A_1$. So $R = A_1 \oplus A_1' \oplus aR = A_1 \oplus A_1'\oplus aA_1 \oplus aY_1$ implying that $$R = K \oplus A_1 \oplus Y_1 = A_1 \oplus (A_1' \oplus aA_1) \oplus aY_1,$$ where $Y_1 \subseteq aR$ and $A_1' \oplus aA_1 \cong A_1'\oplus A_1 = A \cong R/aR$.

\bigskip
Now suppose we have found the right ideals $A_i,\;A_i',\;Y_i$ for $i=1,\dots,n$ such that $(*)$ holds for every $j= 1,\dots, n$ with $Y_i \subseteq a^iR$, $A_i' \oplus aA_i \cong R/aR$ for every $i=1,\dots,n$ and $A_i' \oplus aA_i = A_{i+1}'\oplus A_{i+1}$ for every $i = 1,\dots, n-1$. As $K$ has the exchange property and $$R = K \oplus (\bigoplus_{i=1}^{n}A_i) \oplus Y_n = (\bigoplus_{i=1}^{n}A_i) \oplus (A_n' \oplus aA_n) \oplus aY_n,$$ by Lemma 1 we have decompositions $A_n' \oplus aA_n = A_{n+1} \oplus A_{n+1}'$ and $aY_n = Y_{n+1} \oplus Y_{n+1}'$ such that
$$R = K \oplus (\bigoplus_{i=1}^{n+1}A_i) \oplus Y_{n+1}.$$ So $aR = \bigoplus_{i=1}^{n+1}aA_i \oplus aY_{n+1}$ and $aA_{n+1} \cong A_{n+1}$. Now as $A_j' \oplus aA_j = A_{j+1} \oplus A_{j+1}'$ for every $j = 1,\dots, n$ we have $$R = A_1 \oplus A_1' \oplus aR = A_1 \oplus A_1' \oplus \bigoplus_{i=1}^{n+1} aA_i  \oplus aY_{n+1} = (\bigoplus_{i=1}^{n+1} A_i) \oplus (A_{n+1}' \oplus aA_{n+1}) \oplus aY_{n+1},$$ with $Y_{n+1} \subseteq a^{n+1}R$ and $A_{n+1}' \oplus aA_{n+1} \cong A_{n+1}' \oplus A_{n+1} = A_n' \oplus aA_n \cong R/aR.$ \qed

\bigskip
\noindent{\bf Lemma 3 \cite[Lemma 2.8]{Ni1}.} {\em Let $P_R$ be a projective module and $P = A + B$ where $A \subseteq^{\oplus} P$. Then $B = C\oplus D$ for some submodules $C$ and $D$ such that $P=A\oplus C$.}\\[4mm]
{\bf Theorem 4.} {\em Let $a$ be an element of a ring $R$ such that $a^n$ is regular for every positive integer $n$. Then for every $n \in \mathbb{N}$ there exist decompositions $$ R_R = K \oplus X_n \oplus Y_n = E_n \oplus X_n \oplus aY_n,$$ where $K = \rm{r}_R(a)$, $Y_n \subseteq a^nR$, $K \oplus Y_n = K + a^nR$ and $E_n \cong R/aR$. If  $a^n = 0$, then $K \cong E_n \cong R/aR$ implying that $a$ is unit-regular.}\\[3mm]
{\bf Proof.} For every positive integer $i$ we will inductively construct right ideals $A_i,\;A_i',\;Y_i$ of $R$ such that for every $j\geq 1$, $$R = K \oplus (\bigoplus_{i=1}^{j}A_i) \oplus Y_j = (\bigoplus_{i=1}^{j}A_i) \oplus (A_j' \oplus aA_j) \oplus aY_j, \hspace{10mm} (**)$$ where $Y_j \subseteq a^jR$, $aY_j = a^{j+1}R$ and $A_j' \oplus aA_j = A_{j+1} \oplus A_{j+1}' \cong R/aR$. Then we have the desired decompositions by putting $X_n = \bigoplus_{i=1}^{n}A_i$ and $E_n = A_n' \oplus aA_n$.

\bigskip
Note that the left multiplication by $a$ induces an epimorphism from $R \to aR/a^{n+1}R$ with kernel $K + a^nR$. As $a$ and $a^{n+1}$ are regular, $aR/a^{n+1}R$ is projective implying that $K + a^nR \subseteq^{\oplus} R_R$ for each $n$.  By Lemma 3, $K + aR = K \oplus Y_1$ for some $Y_1 \subseteq aR$. If $Y_1'= aR \cap K$, then on intersecting with $aR$ we have $aR = Y_1' \oplus Y_1$, $a^2R = aY_1$ and $Y_1' \subseteq^{\oplus} R_R$. So $K = Y_1' \oplus A_1'$ for some $A_1'$. Also for some right ideal $A_1$ we have $R = (K+aR) \oplus A_1 = K \oplus Y_1 \oplus A_1 = Y_1' \oplus A_1' \oplus Y_1 \oplus A_1$.  As $K \cap (A_1 \oplus Y_1) = 0$, $aR = aA_1 \oplus aY_1$ and  $aA_1\cong A_1$. So $R = A_1 \oplus A_1' \oplus aR = A_1 \oplus A_1' \oplus aA_1 \oplus aY_1$. Thus $$R = K \oplus A_1 \oplus Y_1 = A_1 \oplus (A_{1}' \oplus aA_1) \oplus aY_1,$$ where $Y_1 \subseteq aR$, $K \oplus Y_1 = K + aR$ and $A_1' \oplus aA_1 \cong A_1' \oplus A_1 \cong R/aR$.

\bigskip
Now suppose we have found the right ideals $A_i,\;A_i',\;Y_i$ for $i=1,\dots,n$ such that we have decompositions as in $(**)$ for every $j= 1,\dots, n$ with $Y_i \subseteq a^iR$, $K \oplus Y_i = K + a^iR$, $A_i' \oplus aA_i \cong R/aR$ for every $i=1,\dots,n$ and $A_i' \oplus aA_i = A_{i+1}'\oplus A_{i+1}$ for every $i = 1,\dots, n-1$. By Lemma 3, $K + a^{n+1}R = K \oplus Y_{n+1}$ for some $Y_{n+1} \subseteq a^{n+1}R$.  Note that $aY_n = a^{n+1}R \subseteq K \oplus Y_{n+1}$ and $(K \oplus Y_{n+1}) \cap \bigoplus_{i=1}^{n}A_i \subseteq (K \oplus Y_n) \cap \bigoplus_{i=1}^{n}A_i = 0$. As $K + a^nR \subseteq^{\oplus} R_R$ for each $n$, $K \oplus Y_{n+1} = K + a^{n+1}R \subseteq^{\oplus} K + a^nR = K \oplus Y_n$ and $K \oplus Y_{n+1}  \oplus \bigoplus_{i=1}^{n}A_i \subseteq^{\oplus} K \oplus Y_n  \oplus \bigoplus_{i=1}^{n}A_i = R.$ Thus $$R = (K \oplus Y_{n+1}) + R = (K \oplus Y_{n+1}) + (\bigoplus_{i=1}^{n}A_i \oplus (A_{n}' \oplus aA_n) \oplus aY_n) = (K \oplus Y_{n+1} \oplus \bigoplus_{i=1}^{n}A_i) + (A_n' \oplus aA_n).$$  Again by Lemma 3, $A_n' \oplus aA_n = A_{n+1} \oplus A_{n+1}'$ such that $ R = K \oplus Y_{n+1} \oplus \bigoplus_{i=1}^{n+1}A_i.$ So $aR = \bigoplus_{i=1}^{n+1}aA_i\oplus aY_{n+1}$, $aA_{n+1} \cong A_{n+1}$.  Now as $A_j' \oplus aA_j = A_{j+1} \oplus A_{j+1}'$ for every $j = 1,\dots, n$ we have $$R = A_1 \oplus A_1' \oplus aR = A_1 \oplus A_1' \oplus \bigoplus_{i=1}^{n+1}aA_i\oplus aY_{n+1} =(\bigoplus_{i=1}^{n+1}A_i) \oplus (A_{n+1}' \oplus aA_{n+1}) \oplus aY_{n+1},$$ with $Y_{n+1} \subseteq a^{n+1}R$, $K \oplus Y_{n+1} = K + a^{n+1}R$ and $A_{n+1}' \oplus aA_{n+1} \cong A_{n+1}' \oplus A_{n+1} = A_n' \oplus aA_n \cong R/aR.$ \qed


  \bigskip
\noindent {\bf Acknowledgements:} We thank Pace P. Nielsen for useful comments that improved the presentation. We also thank the referee for many corrections and improvements.
\bibliographystyle{amsplain}
\bibliography{RegularNilpotentElements}

\end{document}